\documentclass[12pt]{amsart}

\usepackage{amsmath}
\usepackage{amssymb, amsfonts}

\theoremstyle{plain}

\newcommand{\beq}{\begin{equation}}
\newcommand{\eeq}{\end{equation}}

\newtheorem{theor}{Theorem} 
\newtheorem{prop}{Proposition}

\def\S{\sum_{j=0}^7}

\title{The Area Integral and Pointwise Boundary Limits of Monogenic Functions}
\author{  Alexander Kheyfits }

\begin{document}
\maketitle

\begin{abstract}
We prove a criterion of the existence of the non-tangen-tial limit values at a given boundary point for octonion-valued monogenic functions in the half-space. It is also proved that the non-tangential limit for the scalar component of a monogenic function and those limits for all of its other seven components either exist or does not exist simultaneously. The non-associativity and non-commutativity of the octonions can be circumvented in certain problems by making use of P. Stein method employing subharmonic functions.   \end{abstract}

\section{Introduction and Statement of Results}

The classical Fatou theorem \cite{Fa} states that a bounded analytic function in a smooth complex domain has non-tangential boundary values at almost every, with respect to the Lebesgue measure, boundary point. The statement says nothing about the exceptional set of measure zero. It is less known, however, that in the same work Fatou actually gave a sufficient condition for the existence of the normal boundary limit at an individual point; the necessity of that condition was proved by Loomis \cite{Lo} much later. Such statements are now called \emph{pointwise Fatou theorems}.

These results were extended to the harmonic functions in $\mathbb{R}^n, n\geq 2$, caloric functions \cite{Bro}, and functions harmonic with respect to  the stationary Schr\"odinger operators, see \cite{Khe1} and the references therein. In this note, stimulated by \cite{Bro}, we study the boundary properties of the octonion-valued monogenic functions.

\vspace{.5cm}

\footnoterule
\par Key words: \emph{Octonion-valued monogenic functions, Area integral, Non-tangential boundary values, Normal boundary values at a given point, Subharmonicity of powers of monogenic functions}.    \\
\par 2010 Mathematics Subject Classification:  30G35; 31B05; 31B25; 35E99   \\

\pagebreak

Octonionic monogenic functions, due to both  their pure mathematical interest, and to an inherent even though not completely fulfilled utility in physics \cite{Ok,GurTze}, steadily attract attention of researchers. J. Baez formulates in his very informative survey \cite[p. 201]{Ba}  the development of an octonionic analogue of the theory of analytic functions as the first item in his list of 14 open important octonion-related  problems. 

There is a growing literature devoted to the study of the monogenic functions in various settings, see, e.g., \cite{CatKheTep, KheTep, Khe2, LP1, LP2, LP3} and especially \cite{LPQ} and the references therein. In those papers certain important results of classical complex analysis have been extended to the octonion-valued monogenic functions; for instance, the Phragm\'{e}n-Lindel\"{o}f principle,  the three-lines theorem, Paley-Wiener theorem, and an analog of M. Riesz theorem on the boundedness of the conjugation operator, to name just a few.

Here we study the existence of non-tangential boundary values almost everywhere and give a criterion for the existence of finite or infinite boundary limits for these functions at a given boundary point. We concentrate on the function-theoretic properties of the octonion-valued monogenic functions, and leave aside any results about  the Clifford-valued, in particular, quaternionic functions. As for applications to PDEs, see, e.g., the recent monograph \cite{NaTV}. 

However, these generalizations are not trivial, as can be seen, for example, from the fact that unlike the classical case, the octonionic version of the M. Riesz theorem on the conjugate harmonic functions is valid not for all positive $p$, but only for $p\geq 6/7$ -- see \cite{KheTep, CatKheTep}.  \\

To study the non-tangential boundary values of the monogenic functions, we define for these functions an analog of the area integral in the half-space \cite{Ste1}. The area integrals for the monogenic functions in Clifford algebras are known, see, e.g., \cite[Chap. 4]{Mi}, however, there they were connected with the Hardy spaces of Clifford-valued monogenic functions, which are not considered in this note. Now we remind some terminology and state our results, the proofs are given in the next section. 

\pagebreak

The \emph{non-commutative, non-associative, alternative division algebra of octonions} $\mathcal{O}$ is an eight-dimensional vector space with the \emph{basis}  
$\{\mathbf{e}_0\equiv 1, \mathbf{e}_1,\ldots , \mathbf{e}_7 \}  $, that satisfies the multiplication table 

\[\left[ \begin{array}{c|c|c|c|c|c|c|c|c} \hline 
\mathbf{e}_i \times \mathbf{e}_j & \mathbf{e}_0=1  & \mathbf{e}_1 & \mathbf{e}_2 & \mathbf{e}_3 & \mathbf{e}_4 & \mathbf{e}_5 & \mathbf{e}_6 & \mathbf{e}_7 \\
 \hline 
 \hline
 \mathbf{e}_0=1 & 1 & \mathbf{e}_1 & \mathbf{e}_2 & \mathbf{e}_3 & \mathbf{e}_4 & \mathbf{e}_5 & \mathbf{e}_6  &  \mathbf{e}_7   \\ \hline 
\mathbf{e}_1 &  \mathbf{e}_1  & -1 & \mathbf{e}_4 & \mathbf{e}_7 & -\mathbf{e}_2 & \mathbf{e}_6 & -\mathbf{e}_5 & -\mathbf{e}_3 \\ \hline 
\mathbf{e}_2 &  \mathbf{e}_2   & -\mathbf{e}_4  & -1 & \mathbf{e}_5 & \mathbf{e}_1 & -\mathbf{e}_3 & \mathbf{e}_7 & -\mathbf{e}_6 \\ \hline 
\mathbf{e}_3 &  \mathbf{e}_3   & -\mathbf{e}_7 & -\mathbf{e}_5 & -1 & \mathbf{e}_6 & \mathbf{e}_2 & -\mathbf{e}_4 & \mathbf{e}_1 \\ \hline 
\mathbf{e}_4 &  \mathbf{e}_4   & \mathbf{e}_2 & -\mathbf{e}_1 & -\mathbf{e}_6 & -1 & \mathbf{e}_7 & \mathbf{e}_3 & -\mathbf{e}_5 \\ \hline 
\mathbf{e}_5 &  \mathbf{e}_5  & -\mathbf{e}_6 & \mathbf{e}_3 & -\mathbf{e}_2 & -\mathbf{e}_7 & -1 & \mathbf{e}_1 & \mathbf{e}_4 \\ \hline 
\mathbf{e}_6 &  \mathbf{e}_6   & \mathbf{e}_5 & -\mathbf{e}_7 & \mathbf{e}_4 & -\mathbf{e}_3 & -\mathbf{e}_1 & -1 & \mathbf{e}_2 \\ \hline 
\mathbf{e}_7 &  \mathbf{e}_7   & \mathbf{e}_3 & \mathbf{e}_6 & -\mathbf{e}_1 & \mathbf{e}_5 & -\mathbf{e}_4 & -\mathbf{e}_2 & -1 \\ \hline   \end{array}\right]. \vspace{.5cm}\]

In what follows, let $\mathbf{x}=(x_0,x_1,\ldots ,x_7)\in \mathbb{R}^8$ be an eight-dimensional vector of real numbers and
\[\mathbf{o} = \S x_j \mathbf{e}_j \] 
be a generic octonion, $\mathbf{o}\in \mathcal{O}$. Consider eight real-valued continuously differentiable functions 
\[f_0(\mathbf{x}), f_1(\mathbf{x}), \ldots, f_7(\mathbf{x}),\; \mathbf{x}\in \Omega, \]
in a simply-connected domain $\Omega \subset \mathbb{R}^8$. The \emph{octonion-valued left-monogenic functions} in $\Omega$ are defined as eight-dimensional vector-functions 
\beq f(\mathbf{x})=\S  f_j(\mathbf{x}) \mathbf{e}_j, \vspace{.2cm} \mathbf{x}\in \Omega, \eeq
satisfying the operator equation
\beq D[f]=\textbf{0}, \eeq
where 
\[D[f](\mathbf{x})=\S \frac{\partial f_j(\mathbf{x})}{\partial x_j} \mathbf{e}_j  \]
is the Dirac or Cauchy-Riemann-Fueter operator; here $\partial/\partial x_j,\, j=0,1,\ldots ,7, $ are partial derivatives with respect to the coordinates in  $\mathbb{R}^8$. Thus, we study functional-theoretical properties of the elements in the kernel of the Dirac operator $D$.  

Solutions of the system $[f]D=\textbf{0}$ are called right-monogenic functions; the functions, which are both left- and right-monogenic, are called (two-sided) monogenic functions. Hereafter, we always discuss the left-monogenic functions; the proofs hold good for the right- and two-sided monogenic functions as well. Moreover, since adding a constant to any component $f_j$ of any solution of system (3) below gives another solution of the system, we will assume that $f(\mathbf{0})=\mathbf{0}$. \\

The non-commutativity and non-associativity of the algebra of octonions $\mathcal{O}$ result in certain difficulties in a study of this algebra. However, as early as in 1933, P. Stein \cite{StP} employed the subharmonic functions in his study of M. Riesz theorem about the harmonic conjugate functions. Since then, this approach has been used by various authors, see \cite{SteWe, KheTep} and the referenced above. If the final claim depends only on the modulus $|f|$ of the left- or right-monogenic function $f$, then the use of subharmonicity allows in certain problems to fix  an ordering and an association from the outset and work with this order to the end of the proof, when it can be seen that the result does not depend upon a particular ordering and association. 

We systematically use this approach in what follows. \\

Inserting representation (1) into (2) and using the linear independence of the basis octonions $\mathbf{e}_0, \ldots , \mathbf{e}_7$, 
it follows that equation (2) is equivalent to the system of eight first-order linear partial differential equations with constant coefficients with respect to the functions 
$f_0, f_1, \ldots , f_7$. This system  can be written down as the equivalent matrix equation  $\vspace{.5cm}$
\beq \left[\begin{array}{rrrrrrrr} \frac{\partial}{\partial x_0} &
-\frac{\partial}{\partial x_1} & -\frac{\partial}{\partial x_2} &
-\frac{\partial}{\partial x_3} & -\frac{\partial}{\partial x_4} &
-\frac{\partial}{\partial x_5} & -\frac{\partial}{\partial x_6} & -\frac{\partial}{\partial x_7} \vspace{.24cm}\\

\frac{\partial}{\partial x_1} & \frac{\partial}{\partial x_0} &
-\frac{\partial}{\partial x_4} & -\frac{\partial}{\partial x_7} &
\frac{\partial}{\partial x_2} & -\frac{\partial}{\partial x_6} &
\frac{\partial}{\partial x_5} & \frac{\partial}{\partial x_3} \vspace{.24cm}\\

\frac{\partial}{\partial x_2} & \frac{\partial}{\partial x_4} &
\frac{\partial}{\partial x_0} & -\frac{\partial}{\partial x_5} &
-\frac{\partial}{\partial x_1} & \frac{\partial}{\partial x_3} &
-\frac{\partial}{\partial x_7} & \frac{\partial}{\partial x_6} \vspace{.24cm}\\

\frac{\partial}{\partial x_3} & \frac{\partial}{\partial x_7} &
\frac{\partial}{\partial x_5} & \frac{\partial}{\partial x_0} &
-\frac{\partial}{\partial x_6} & -\frac{\partial}{\partial x_2} &
\frac{\partial}{\partial x_4} & -\frac{\partial}{\partial x_1} \vspace{.24cm}\\

\frac{\partial}{\partial x_4} & -\frac{\partial}{\partial x_2} &
\frac{\partial}{\partial x_1} & \frac{\partial}{\partial x_6} &
\frac{\partial}{\partial x_0} & -\frac{\partial}{\partial x_7} &
-\frac{\partial}{\partial x_3} & \frac{\partial}{\partial x_5} \vspace{.24cm}\\

\frac{\partial}{\partial x_5} & \frac{\partial}{\partial x_6} &
-\frac{\partial}{\partial x_3} & \frac{\partial}{\partial x_2} &
\frac{\partial}{\partial x_7} & \frac{\partial}{\partial x_0} &
-\frac{\partial}{\partial x_1} & -\frac{\partial}{\partial x_4} \vspace{.24cm}\\

\frac{\partial}{\partial x_6} & -\frac{\partial}{\partial x_5} &
\frac{\partial}{\partial x_7} & -\frac{\partial}{\partial x_4} &
\frac{\partial}{\partial x_3} & \frac{\partial}{\partial x_1} &
\frac{\partial}{\partial x_0} & -\frac{\partial}{\partial x_2} \vspace{.24cm}\\

\frac{\partial}{\partial x_7} & -\frac{\partial}{\partial x_3} &
-\frac{\partial}{\partial x_6} & \frac{\partial}{\partial x_1} &
-\frac{\partial}{\partial x_5} & \frac{\partial}{\partial x_4} &
\frac{\partial}{\partial x_2} & \frac{\partial}{\partial x_0}
\end{array} \right] \times
\left[\begin{array}{c}
  f_0 \vspace{.3cm} \\
  f_1 \vspace{.3cm} \\
  f_2 \vspace{.3cm} \\
  f_3 \vspace{.3cm} \\
  f_4 \vspace{.3cm} \\
  f_5 \vspace{.3cm} \\
  f_6 \vspace{.3cm} \\
  f_7
\end{array}  \right]
=\textbf{0}. \vspace{.5cm} \eeq with respect to the vector-function $f(\mathbf {x})$.

Differentiating the equations of system (3) and combining them in an obvious way, one derives the eight equations
\[\Delta f_0(\mathbf{x})=\cdots =\Delta f_7(\mathbf{x}) =0,\; \mathbf{x}\in \Omega, \]
where $\Delta$ is the eight-dimensional Laplace operator. Hence, as it is well-known, all the components $f_0,\ldots ,f_7$ of an octonion-valued monogenic function $f$ are the classical harmonic functions.

Equations (2) such that each component $f_j,\; j\geq 0$, is harmonic, are called the Generalized Cauchy-Riemann systems (GCR) - see Stein and Weiss \cite[pp. 231-234]{SteWe}. More general systems 
\[ \sum_{j=0}^n A_j \frac{\partial f}{\partial x_j}+Bf=\textbf{0} \]
with constant matrices $A_j, \; j=0,1,\ldots ,n,$ and $B$ were considered, in different context, by Evgrafov \cite{Evg}.

Following the insight of Calder\'on, Stein and Weiss have proved that for any GCR system $F$ there exists a nonnegative index 
$p_0<1$ such that $|F|^p$ is a subharmonic function for all $p\geq p_0$. In particular, it is known 
\cite[p. 234]{SteWe} that for the M. Riesz system in $\mathbb{R}^n \vspace{.2cm}$,
\[\begin{cases} \frac{\partial f_1}{\partial x_1} + \cdots + \frac{\partial f_n}{\partial x_n} = 0   
\vspace{.3cm}  \\
\frac{\partial f_i}{\partial x_j}=\frac{\partial f_j}{\partial x_i},\; i,j=1,\ldots ,n, \end{cases} \] 
the exact value is $p_0=(n-2)/(n-1)$. Due to the non-associativity of the octonions, this statement is not immediately obvious for the latter. However, the same assertion is valid for system (2)-(3) in $\mathbf{R}^8$ as well; namely, it has been proven in \cite{KheTep} that for the octonion-valued monogenic functions, that is, for the solutions of system (2)-(3),
\[p_0=\frac{n-2}{n-1}\biggl|_{n=8}=6/7.\]

As Stein and Weiss have noticed (ibid., p. 233) the inequality $p_0<1$ allows one to develop a substantive theory of the Hardy spaces for the corresponding systems (2), in our case for octonionic monogenic functions. This inequality is also important for many other problems, see, e.g., applications in \cite{CatKheTep}. 

We consider monogenic functions on upper half-space 
\[\mathbb{R}^8_+ = \Bigg\{\mathbf{x}=(x_0,X): \; x_0>0,\; X=( x_1,\ldots ,x_7)\in \mathbb{R}^7  \Bigg\}. \]
For a fixed point $Y \in \mathbb{R}^7$ and positive numbers $h>0,\; \alpha >0,$ define an (open) truncated cone \cite[p. 138]{Ste1}
\[\Gamma (Y)=\{\mathbf{x}=(x_0,X): \; |Y - X| < \alpha x_0,\; 0<x_0<h.  \]

Remind the following definition by E. Stein \cite[p. 143]{Ste1}:  \\

The function $f(\mathbf{x}) = (f_0(\mathbf{x}), f_1(\mathbf{x}),\ldots ,f_7(\mathbf{x}))$ defined in the half-space $\mathbb{R}^8_+$, has a non-tangential limit $l=(l_0,l_1,\ldots ,l_7)$ at the point $(0,Y)\in \partial \mathbb{R}^8_+$, if $f(\mathbf{x})\rightarrow l$ for every fixed $\alpha$, that is, $f_j(\mathbf{x})\rightarrow l_j $ for $j=0,1,\ldots ,7$, as $x_0\rightarrow 0$ and $ |Y-X|<\alpha x_0$. \\

In the same manner, we call the function $f=(f_0, f_1,\ldots , f_7)$ defined in the half-space $\mathbb{R}^8_+$, \emph{non-negative}, if each its component $f_j,\; j=0,1,\ldots ,7$, is \emph{non-negative} in the half-space $\mathbb{R}^8_+$, $f_j \geq 0$. The following 64-tuple of the partial derivatives of the components $f_0,\ldots ,f_7$ is called the \emph{second gradient} of $f$, Cf. \cite [p.162]{Ste1}, \[\nabla f(\mathbf{x}) = (\partial f_0/\partial x_0, \partial f_0/\partial x_1, \ldots , \partial f_0/\partial x_7, \partial f_1/\partial x_0, \ldots , \partial f_7/\partial x_7); \]
their arrangement is immaterial, since we only use the magnitude 
\beq |\nabla f| ^2= \sum_{j=0}^7 \sum_{k=0}^7 \biggl|\frac{\partial f_j}{\partial x_k} \biggl|^2.  \eeq

To state our results about the boundary limits of the octonionic monogenic functions, we define now a "monogenic" version of the classical \emph{area integral}.  \\

Given an octonion-valued monogenic function $f(\mathbf{x})$ in the half-space $\mathbb{R}^8_+$, its area integral or the Lusin square function is
\beq  \left(\mathcal{A}_f(X)\right) ^2=  \int \int_{\Gamma (X)} y_0^{-6} |\nabla f((\mathbf{y})|^2 d\mathbf{y},\;\; \mathbf{x} =(x_0,X)\in \mathbb{R}^8_+. \eeq
This is a real-valued integral of positive functions, hence, commutativity-associativity-related issues do not occur here. \\

Now we can state the octonionic analogs of the results \cite{Ste1} connecting the area integral with boundary properties of the harmonic functions; see \cite{Bro} for "caloric" version of these results in a Lipschitz domain.
\begin{prop} Let $f$ be a (left) octonion-valued monogenic function in the half-space $\mathbb{R}^8_+$ and $E\subset \mathbb{R}^7 = \partial \mathbb{R}^8_+$ be a bounded set at the boundary of the half-space.

1. If the area integral $\mathcal{A}_f(Y)$ is finite for every point $Y\in E$, then function $f$ is locally bounded in a neighborhood of $E$ and has a non-tangential boundary limit at almost every point in $E$.

2. On the other hand, if for every boundary point $Y\in E$, $f$ is bounded in a cone $\Gamma (Y)$, then the area integral $ \mathcal{A}_f (Y) $ is finite for almost every $Y\in E$.
\end{prop}  

An octonion-valued monogenic function is a vector-function, however, as the next theorem shows, in certain problems, the existence of the boundary limit of scalar component is equivalent to the simultaneous limits of the other components. \\

\begin{theor} Under the same assumptions as in Proposition 1, if its scalar component $f_0$ has non-tangential limits at every point of the boundary set $E$, then all the functions $f_1,f_2, \ldots ,f_7$ have the non-tangential boundary limits at a.e. point of $E$.

Conversely, if all the functions $f_1,f_2, \ldots ,f_7$ have simultaneously the non-tangential boundary limits 
at every point of $E$, then $f_0$ has non-tangential limits at a.e. points of the boundary set $E$. 
\end{theor}

In the classical theory of harmonic functions an important role is played by the notion of  the conjugate harmonic functions. For a monogenic function $f$, the function $f_0$ can be considered a given harmonic function -- a \emph{scalar} component of $f$, while the vector-function $f_v=(f_1,\ldots ,f_7)$ is the \emph{vector} component, a conjugate function. As Stein noticed \cite[p. 163]{Ste1}, the theorem cannot be improved in that no part of the vector component can be omitted from the assumptions of the second statement of the theorem.  \\

To state our next result, we remind that by the Riesz-Herglotz theorem, a positive harmonic function $u$ 
in the half-space $\mathbb{R}^n_+$ has an integral  representation
\beq u(\mathbf{x})=u(x_0, X)=cx_0+ \int_{\mathbb{R}^{n-1}} P(x_0,X-Y) d\mu(Y), \eeq
where $c\geq 0$ is a non-negative constant, $P(x_0,X)$ the Poisson kernel for the half-space 
$\mathbb{R}^n_+$, and $d\mu$ is the boundary measure of the function $u$, i.e., a non-negative Radon measure on $\mathbb{R}^{n-1}$. In the case of a half-plane $\mathbb{R}^2_+$, this measure is a non-decreasing function, and a non-negative harmonic function in a half-plane has a normal boundary limit at a boundary point $X$ if and only if the function $\mu(t)$ has the \emph{symmetric derivative} at the point $X$, i.e., 
\[\lim_{t\rightarrow 0^+} (2t)^{-1}\{\mu(X+t)-\mu(X-t)\}. \]

In the $n-$dimensional case, the criterion is given in terms of the total Riesz mass of the boundary ball $B(X,t)\subset \mathbb{R}^{n-1} $ centered at the point $X\in \partial \mathbb{R}^{n-1}$, i.e.
\[n(t, X)=\int_{B(X,t)} d\mu(Y). \]

\textbf{Theorem K} (\cite{Khe1} and the references therein). \emph{Let $u$ be a positive harmonic function in the half-space $\mathbb{R}^n_+$. In order for the finite limit 
\[\lim_{x_0\rightarrow 0^+} u(x_0,X)=l,\;  0\leq l < \infty,   \] 
to exist at a given boundary point $X$, it is necessary and sufficient that 
\[\lim_{t\rightarrow 0^+} t^{1-n} n(t,X)= \omega_{n-1} l,  \] 
where $\omega_{n-1}$ is the volume of the unit ball in $\mathbb{R}^{n-1}$. }  \\

\emph{ In the case $l=+\infty$, the necessary and sufficient condition is the existence of the limit (for any, and therefore, for every $a>0$)
\[ \lim_{t \rightarrow 0^+} t \int_t ^a  t^{-n-1} n(t,X) dt = +\infty.  \]  }

The last result of this note generalizes this assertion onto the functions under consideration. Let 
$f(\mathbf{x})=(f_0, f_1,\ldots ,f_7)$ be an octonionic monogenic non-negative function in the half-space 
$\mathbb{R}^8_+$, thus $f_j(\mathbf{x})\geq 0,\; j=0,1,\ldots, 7$. Hence each component 
$f_j,\; 0\leq j\leq 7,$ has the Riesz-Herglotz integral representation similar to (6) with the measure 
$d\mu_j(Y)$ satisfying Theorem K above. By the Riesz associated measure of the monogenic function $f(\mathbf{x})$ we understand the vector-measure 
\[d\mu_f(Y)= (d\mu_0(Y),\ldots, d\mu_7(Y)).\]
The associated mass of the function $f(\mathbf{x})$ in a ball is defined similarly, as the vector of the integrals
\[n_f(t, X)= (n_0(t,X), n_1(t,X),\ldots , n_7(t,X))= \] 
\[= \int_{B(X,t)} d\mu_f(Y)= \]
\[=\bigg(\int_{B(X,t)} d\mu_0(Y),\ldots , \int_{B(X,t)}d\mu_7(Y) \bigg). \]

\begin{theor} Let $f$ be a positive octonion-valued monogenic function in the upper half-space $\mathbb{R}^8_+$, and $l=(l_0,l_1,\ldots ,l_7) $ be a vector with non-negative, maybe infinite components. In order for the limit 
\[\lim_{x_0\rightarrow 0^+} f(\mathbf{x})=l   \] 
to exist at a given boundary point $(0, X),\; \mathbf{x}=(x_0,X)$, 
it is necessary and sufficient that for every $j=0,1,\ldots ,7$, either 
\[\lim_{t\rightarrow 0^+} t^{-7} n_j(t,X)= \omega_7 l_j   \] 
if the corresponding $l_j$ is finite, or
\[\lim_{t \rightarrow 0^+} t \int_t ^a  t^{-9} n_j(t,X) dt = +\infty  \] 
if $l_j= +\infty$; here $\omega_7$ is the volume of the unit ball in $\mathbb{R}^7$.  \end{theor}

\section{Proofs} 

\emph{Proof of Proposition 1}. Definitions (4) and (5) imply the equation
\[\mathcal{A}_f(Y)^2 = \mathcal{A}_{f_0}(Y)^2+\mathcal{A}_{f_1}(Y)^2+\cdots + \mathcal{A}_{f_7}(Y)^2. \]
Thus, if the area integral $\mathcal{A}_f(Y)$ is finite, then all the partial area integrals $\mathcal{A}_{f_j}(Y)$ are finite, whence by the theorem of E. Stein \cite[Theorem 1]{Ste1}, every (harmonic) component of the function $f$ is upper-bounded and has non-tangential boundary limits at almost every boundary point. A similar reasoning proves the converse statement of the proposition.      
$\hfill  \qed$ \\

\emph{Proof of Theorem 1.} In the proof we essentially use the following fundamental theorem by E. Stein \cite[p. 161, Theorem 3]{Ste1}.  \\

\textbf{Theorem S}. \emph{Let $u$ and $v$ be vectors of harmonic functions in $\mathbb{R}^{n+1}_+$ of 
$k$ and $m$ components respectively, connected by the equation 
\beq \frac{\partial^r u}{\partial x^r_0}=P(D) v, \eeq
where $P(D)$ is a $k\times m$ matrix whose entries are differential polynomials (with constant coefficients) homogeneous of degree $r$, $r\geq 1$. Suppose that for a given set $E$, $E\subset \mathbb{R}^n$, $v$ has a non-tangential limit for every $X \in E$. Then $u$ has a non-tangential limit for a.e. $X \in E$. } \\

Let $u=(f_0)$ be a one-by-one matrix function, 
\[v=\left( \begin{array}{c}f_1 \\
f_2 \\ \vdots  \\ f_7 \end{array} \right) \]
a seven-by-one column matrix, and 
\[P(D)=\left( \frac{\partial}{\partial x_1}, \ldots , \frac{\partial}{\partial x_7} \right) \]
a vectorial, one-by-seven, linear differential operator of order $r=1$ with constant coefficients. 

In these notations, the first equation of system (3),
\[\frac{\partial f_0}{\partial x_0} = \sum_{j=1}^7 \frac{\partial f_j}{\partial x_j}, \]
can be written as
\[\frac{\partial u}{\partial x_0}=P(D)v. \]
The second part of our Theorem 1 now immediately follows from Theorem S above.  \\

The first part of Theorem 1 also follows from Theorem S, but instead of the first equation of system (3), we use its second, third, etc., equations, and must choose another operator $P(D)$. The second equation of system (3) implies the equation
\[\frac{\partial f_0}{\partial x_1} = - \frac{\partial f_1}{\partial x_0}+ \frac{\partial f_2}{\partial x_4} + 
\frac{\partial f_3}{\partial x_7}-\frac{\partial f_4}{\partial x_2}+\frac{\partial f_5}{\partial x_6} - 
\frac{\partial f_6}{\partial x_5} - \frac{\partial f_7}{\partial x_3}. \]
The third equation of (3) gives 
\[\frac{\partial f_0}{\partial x_2} = - \frac{\partial f_1}{\partial x_4} - \frac{\partial f_2}{\partial x_0}  + 
\frac{\partial f_3}{\partial x_5} + \frac{\partial f_4}{\partial x_1} - \frac{\partial f_5}{\partial x_3} + 
\frac{\partial f_6}{\partial x_7} - \frac{\partial f_7}{\partial x_6}. \]
The other equations of (3) give similar expressions for $\frac{\partial f_0}{\partial x_j},\; 3\leq j\leq 7$. 

Now if $u=(f_0)$ is a one-by-one matrix-function, then the left-hand side of these equations form a seven-by-one column-vector 
\[U=\left( \begin{array}{c}\frac{\partial f_0}{\partial x_1}  \vspace{.2cm} \\
\frac{\partial f_0}{\partial x_2}   \vspace{.2cm} \\ \vdots   \vspace{.2cm} \\ \frac{\partial f_0}{\partial x_7}  \end{array} \right). \]
Similarly, if $V$ is a seven-by-one column-vector
\[V=\left( \begin{array}{c}f_1   \vspace{.1cm} \\
f_2    \vspace{.1cm} \\ \vdots    \vspace{.1cm} \\ f_7 \end{array} \right), \]
then the system above for the derivatives $\frac{\partial f_0}{\partial x_j},\; 1\leq j\leq 7,$ can be written as a matrix equation
 \[U=P(D)V,\]
where the differential operator $P(D)  \vspace{.2cm}$ is given by the matrix
\[ \left[\begin{array}{rrrrrrr} 
- \frac{\partial}{\partial x_0} &
\frac{\partial}{\partial x_4} & \frac{\partial}{\partial x_7} &
-\frac{\partial}{\partial x_2} & \frac{\partial}{\partial x_6} &
-\frac{\partial}{\partial x_5} & -\frac{\partial}{\partial x_3} \vspace{.24cm}\\

- \frac{\partial}{\partial x_4} &
-\frac{\partial}{\partial x_0} & \frac{\partial}{\partial x_5} &
 \frac{\partial}{\partial x_1} & -\frac{\partial}{\partial x_3} &
\frac{\partial}{\partial x_7} & -\frac{\partial}{\partial x_6} \vspace{.24cm}\\

- \frac{\partial}{\partial x_7} &
-\frac{\partial}{\partial x_5} & -\frac{\partial}{\partial x_0} &
\frac{\partial}{\partial x_6} & \frac{\partial}{\partial x_2} &
-\frac{\partial}{\partial x_4} & \frac{\partial}{\partial x_1} \vspace{.24cm}\\

 \frac{\partial}{\partial x_2} &
-\frac{\partial}{\partial x_1} & -\frac{\partial}{\partial x_6} &
 -\frac{\partial}{\partial x_0} & \frac{\partial}{\partial x_7} &
 \frac{\partial}{\partial x_3} & -\frac{\partial}{\partial x_5} \vspace{.24cm}\\

- \frac{\partial}{\partial x_6} &
\frac{\partial}{\partial x_3} & -\frac{\partial}{\partial x_2} &
-\frac{\partial}{\partial x_7} & -\frac{\partial}{\partial x_0} &
\frac{\partial}{\partial x_1} & \frac{\partial}{\partial x_4} \vspace{.24cm}\\

 \frac{\partial}{\partial x_5} &
-\frac{\partial}{\partial x_7} & \frac{\partial}{\partial x_4} &
-\frac{\partial}{\partial x_3} & -\frac{\partial}{\partial x_1} &
-\frac{\partial}{\partial x_0} & \frac{\partial}{\partial x_2} \vspace{.24cm}\\

 \frac{\partial}{\partial x_3} &
\frac{\partial}{\partial x_6} & -\frac{\partial}{\partial x_1} &
\frac{\partial}{\partial x_5} & -\frac{\partial}{\partial x_4} &
-\frac{\partial}{\partial x_2} & -\frac{\partial}{\partial x_0}
\end{array}     \right].   \vspace{.2cm}  \]

The first part of Theorem 1 now follows from Theorem S.    $\hfill \qed $   \\

\emph{Proof of Theorem 2} follows from Theorem K above similarly to the proof of Proposition 1. Indeed, if the function $f$ has a normal limit $l$ at a boundary point, then by definition, all the components $f_j$ have the normal limits $l_j$, finite or infinite. Therefore, the boundary measures have the corresponding limits by Theorem K. Since this reasoning is reversible, the conclusion follows.       $\hfill \qed$    
\bigskip

\textbf{Acknowledgement}. This work was partially supported by the City University of New York.


$\vspace{.1cm}$

The Graduate Center and Bronx Community College of 

The City University of New York 

Alexander.Kheyfits@gmail.com  

\end{document}